\documentclass[11pt]{article}

\usepackage{pb-diagram}
\usepackage{epsfig}
\usepackage{amsmath}
\usepackage{amssymb}
\usepackage{amsthm}

\theoremstyle{plain}
\newtheorem{theorem}{Theorem}[section]
\newtheorem{lemma}[theorem]{Lemma}
\newtheorem{proposition}[theorem]{Proposition}

\theoremstyle{definition}

\newtheorem{remark}[theorem]{Remark}

\newcommand\bC{{\mathbb C}}

\newcommand\bR{{\mathbb R}}

\newcommand\bZ{{\mathbb Z}}

\renewcommand{\Box}{\blacksquare}

\newcommand{\A}{\mathcal A}

\newcommand{\f}{\mathcal F}

\newcommand{\p}{{\mathcal P}}

\newcommand{\R}{\mathcal R}
\newcommand{\s}{{\mathcal S}}

\newcommand{\ra}{\rightarrow}
\newcommand{\hra}{\hookrightarrow}

\def\inpr{\mathbin{\hbox to 6pt{\vrule height0.4pt width5pt depth0pt \kern-.4pt \vrule height6pt width0.4pt depth0pt\hss}}}

\newcommand{\bc}{{\bf c}}

\newcommand{\ii}{{\bf i}}

\newcommand{\ve}{{\varepsilon}}
\newcommand{\vfi}{{\varphi}}

\begin{document}

\title{On the space of Fredholm operators}

\author{Liviu I. Nicolaescu\\University of Notre Dame\\Notre Dame, IN 46556\\nicolaescu.1@nd.edu}

\date{}



\maketitle

\begin{abstract} We compare various topologies on the space of   (possibly unbounded) Fredholm  selfadjoint operators and explain their  $K$-theoretic relevance.
\end{abstract}


\section*{Introduction}
\setcounter{equation}{0}

The  work of Atiyah and Singer on the index of elliptic operators on  manifolds     has  singled out the role of the space of bounded Fredholm operators  in topology. It is a classifying space for a very useful functor, the topological $K$-theory. This means that  a  continuous family $(L_x)_{x\in X}$ of elliptic  pseudo differential operators  parameterized by a compact $CW$-complex  $X$ naturally defines an element  in the group $K(X)$, the index of the family.

In  most examples, the elliptic operators     are  not bounded operators   and thus  the notion of  continuity has to be  defined  carefully. The operator   theorists have    come up with a   quick  fix. The family $x\mapsto L_x$ of Fredholm operators is  called {\em Riesz  continuous} if and only if the families of {\em bounded operators}
\[
x\mapsto L_x(1+L^*_xL_x)^{-1/2},\;\;x\mapsto L^*_x(1+ L_xL^*_x)^{-1/2}
\]
are continuous with respect to the operator norm.    In concrete applications    this approach can be a nuisance.  For  example, consider as in \cite{N} a  {\em Floer family}  of  {\em elliptic} boundary   value problems (parameterized  by $s\in S^1$)
\begin{equation*}
u(t): [0,1]\ra  {\bC},\;s\in [0,2\pi]\;\;\left\{
\begin{array}{ccc} \;\; \frac{du(t)}{dt}+ a(t) u(t)= 0 &{\rm if}& t\in (0,1)\\
& & \\
u(0)\in {\bR} &  {\rm and} & e^{\ii s}u(1)\in {\bR}
\end{array}
\right.
\tag{$BV_s$}
\label{tag: bvs}
\end{equation*}
where $a:[0,1]\ra {\bC}$ is a given smooth function.  This family ought to be considered continuous but  verifying the above definition can be  quite demanding.   The  first technical goal  of this paper is to elucidate this continuity issue.

As observed   in \cite{AS, Kar}, for $K$-theoretic purposes it suffices to investigate only  (possibly ${\bZ}_2$-graded) selfadjoint operators  (super-)commuting with some Clifford algebra action.    For example, the space of Fredholm  operators on a Hilbert space $H$ can be identified with the space of  odd, selfadjoint Fredholm operators on the ${\bZ}_2$-graded space  $H\oplus H$  via the correspondence
\[
L \mapsto\left[
 \begin{array}{cc} 0 & L^*\\
L & 0
\end{array}
\right].
\]
That is why  we will focus  exclusively on  selfadjoint operators.

In \cite{N} we have argued that in many instances  it is much more
convenient to look at the graphs of Fredholm selfadjoint operators
on a Hilbert space $H$.    If $T$ is such an operator and
$\Gamma_T\subset H\oplus H$ is its graph,  then $\Gamma_T$ is a
Lagrangian  subspace of $H\oplus H$ (with respect to a natural
symplectic structure)  and moreover, the pair $(H\oplus 0,\Gamma_T)$ is Fredholm.   As shown in \cite{N0}, the space of
Fredholm pairs     of Lagrangian subspaces is  a classifying
space for $KO^1$.   (A similar description is valid for all the
functors $KO^n$; see \cite{N}.)

A natural question arises.   Suppose that  two
{\em families of subspaces} determined by the graphs  of two families
of Fredholm operators  are homotopic inside the larger space of Fredholm
pairs of Lagrangian subspaces. Can we conclude that the corresponding
families of Fredholm {\em operators} are also homotopic inside the
smaller space of operators?

The is the second issue we want  to address in this paper.  We will consider  various topologies on the  space of closed,  {\em  unbounded} Fredholm operators  and analyze when  the above graph  map $T\mapsto \Gamma_T$   from operators to  subspaces is a {\em homotopy equivalence}.   Surprisingly, to answer this question we only need to decide the continuity of  Floer type families of boundary value problems. The symplectic reduction technique developed in \cite{N}  coupled with the Bott periodicity will take care of the rest.

The paper consists of three sections.  In Section \ref{s: 1} we   compare two topologies  on the space of  unbounded  Fredholm operators: the gap  topology, given by the gap distance between the graphs, and the Riesz topology,  described above. In the second section  we prove  a general  criterion (Proposition \ref{prop: bvp}) for recognizing  when a   family of boundary value problems, such as (\ref{tag: bvs}), is continuous with respect to the Riesz topology.  In the last section we  address  the connections with $K$-theory.

\bigskip

\noindent {\bf Acknowledgments} This paper addresses  some subtle  omissions in \cite{N}.  I am grateful to Bernhelm Booss-Bavnbek for    his warm reception at Roskilde  University,  where the ideas in this paper were born, and  for  the lively discussions concerning the results of \cite{N}.

\section{Topologies on the space of selfadjoint operators}
\setcounter{equation}{0}
\label{s: 1}

Let $H$ be a separable real Hilbert space.  Denote by $\s$ the space of densely defined, selfadjoint operators on $H$ and by $\mathcal{BS}$ the space  of bounded selfadjoint operators $T:H\ra H$.  Set
\[
[\mathcal{BS}]:=\{T\in\mathcal{BS};\;\;  \|T\|< 1\}.
\]

The {\em Riesz map} is the bijection
\[
\Psi: {\s}\ra[\mathcal{BS}],\;\;A\mapsto A(1+A^2)^{-1/2}.
\]
There are two natural metrics on ${\s}$:   the {\em gap metric}
\[
\gamma(A_0, A_1):= \|(\ii +A_0)^{-1}-(\ii+A_1)^{-1}\|+\|(\ii -A_0)^{-1}-(\ii -A_1)^{-1}\|,
\]
and the {\em Riesz metric}
\[
\rho(A_0, A_1):=\|\Psi(A_0)-\Psi(A_1)\|.
\]

\begin{remark} According to \cite[Thm. IV.2.23]{Kato} we have $\gamma(A_n, A)\ra 0$ if and only if
\[
\delta(\Gamma_{A_n}, \Gamma_A)\ra 0
\]
where $\Gamma_T$ denotes the graph of the linear operator $T$ and $\delta$ denotes the gap between two closed subspaces.
\end{remark}

\begin{lemma}  The identity map $(\s, \rho)\ra (\s, \gamma)$ is continuous.
\label{lemma: ga-ri}
\end{lemma}

\noindent{\bf Proof}\hspace{.3cm}  Observe that for every $A\in \s$ we have
\[
\frac{1}{\ii \pm A}= \frac{A\mp\ii}{1+A^2}= \frac{A}{1+A^2}\mp \frac{1}{1+A^2}=\frac{1}{(1+A^2)^{1/2}}\Psi(A)\mp\ii\frac{1}{1+A^2}
\]
and
\[
\frac{1}{1+A^2}= 1-\Psi(A)^2
\]
so that $\|\Psi(A_n)-\Psi(A)\|\ra 0$ implies $\|(\ii \pm A_n)^{-1}-(\ii \pm A)^{-1}\|\ra 0$.  $\Box$

\bigskip

Denote by $\A$  the  $C^*$-algebra of continuous functions $f: {\bR}\ra {\bC}$   such that the limits
\[
f(\pm\infty):=\lim_{\lambda\ra \pm \infty} f(\lambda)\in {\bC}
\]
exist. Denote by $\A_0$ the subalgebra defined by the  condition
\[
f\in \A_0\Longleftrightarrow f(-\infty)=f(\infty).
\]
Define $P_0, P_\pm \in \A_0$ by
\[
P_0(\lambda)\equiv 1,\;\;P_\pm(\lambda)= (\lambda \pm \ii)^{-1}.
\]
The Stone-Weierstrass approximation theorem shows that the algebra  $\p$ generated by $P_0, P_\pm$ is dense in $\A_0$.

 The functional calculus for selfadjoint operators show that any $A\in \s$ defines a continuous  morphism of $C^*$-algebras
 \[
 \A\ra \mathcal{BS},\;\;f\mapsto f(A).
 \]
\begin{proposition} The following statements are equivalent.

\noindent (i) $\gamma(A_n, A)\ra 0$.

\noindent (ii) $\|f(A_n)-f(A)\|\ra 0$, $\forall f\in \A_0$.
\label{prop: gap}
\end{proposition}

\noindent {\bf Proof}\hspace{.3cm}  Clearly (ii) $\Longrightarrow$ (i) since $P_\pm \in \A_0$ and
\[
\gamma(A_n, A)=\|P_-(A_n)-P_-(A)\|+\|P_+(A_n)-P_+(A)\|.
\]
To prove  (i) $\Longrightarrow$ (ii) we use an idea in \cite[Chap. VIII]{RS}.  Clearly if $\gamma(A_n, A)\ra 0$ then
\[
\|P(A_n)-P(A)\|\ra 0,\;\;\forall P\in \p.
\]
Fix  $f\in \A_0$. Since $\p$ is dense in $\A_0$, for  every $\ve >0$ we can find $P\in \p$ such that $\|f-P\|\leq \ve/3$  and then $n(\ve)>0$ such that, $\forall n\geq n(\ve)$ such that
\[
\|P(A_n)-P(A)\|\leq \ve /3.
\]
Then, $\forall n\geq n(\ve)$ we have
\[
\|f(A_n)-f(A)\|\leq \|f(A_n)-P(A_n)\|+\|P(A_n)-P(A)\|+\|P(A)-f(A)\|\leq \ve.  \;\;\Box
\]

\bigskip

\begin{proposition} Fix a function $\alpha \in \A$ such that $\alpha (\lambda)\equiv 1$ for $\lambda \gg 1$ and $\alpha (\lambda )\equiv 0$ if $\lambda \ll -1$. Then the following  statements are equivalent.

\noindent (i) $\rho(A_n, A)\ra 0$

\noindent (ii) $\|f(A_n)-f(A)\|\ra 0$,  $\forall f\in \A$.

\noindent (iii) $\gamma(A_n, A)\ra 0$ and $\|\alpha(A_n)-\alpha(A)\|\ra 0$.

\label{prop: riesz}
\end{proposition}

\noindent {\bf Proof}\hspace{.3cm} Define $r\in \A$ by
\[
r(\lambda):=\frac{\lambda}{(1+\lambda^2)^{1/2}}.
\]
The equivalence (i) $\Longleftrightarrow$ (ii)  follows exactly as in the proof  of Proposition \ref{prop: gap} using Lemma \ref{lemma: ga-ri} and  the fact that the subalgebra spanned by $\A_0$ and $r$ is dense in $A$. The equivalence (ii) $\Longleftrightarrow$ (iii)  relies on Proposition \ref{prop: gap}  and the fact that the algebra spanned by $\A_0$ and $\alpha$ is dense in $\A$. $\Box$

\bigskip

\begin{remark}{\bf (B. Fuglede)} The  topological spaces $(\s, \rho)$ and $(\s, \gamma)$ are not homeomorphic.   Using Proposition \ref{prop: riesz} it is easy to construct an example  of a sequence $A_n\stackrel{\gamma}{\ra} A$   such that $A_n$ does not  converge to $A$ in the Riesz metric. More precisely consider the space
\[
\ell^2=\Bigl\{ (x_j)_{n\geq 1};\;\;x_j\in {\bR},\;\sum_j x_j^2 <\infty\Bigr\}
\]
with canonical Hilbert basis ${\bf e}_1, {\bf e}_2, \cdots$. For $n=0,1,2,\cdots$ define
\[
A_n: D(A_n)\subset \ell^2\ra \ell^2, \;\;D(A_n)=\Bigl\{(x_j)_{j\geq 1}\in \ell^2;\;\;\sum_{j\geq 1}j^2 |x_j|^2 < \infty\Bigr\}
\]
\[
A_n{\bf e}_j=\left\{
\begin{array}{lr}
j{\bf e}_j, & j\neq n\\
-n{\bf e}_j, & j=n
\end{array}
\right.
\]
One can see that
\[
\|(\ii \pm A_n)^{-1}-(\ii\pm A_0)^{-1}\|= \Bigl|\frac{1}{\ii +n}-\frac{1}{\ii -n}\Bigr|\ra 0
\]
so that $\gamma(A_n, A_0)\ra 0$. On the other hand,  if $\alpha\in \A$ is as in  Proposition \ref{prop: riesz} then for all sufficiently large $n$ we have
\[
\|\alpha(A_n)-\alpha(A_0)\|=1.
\]
\end{remark}

We  now want to present a simple criterion  of $\rho$-convergence. For any closed densely defined operator we denote by $\R(T)\subset {\bC}$ its resolvent set.

\begin{proposition}  Suppose $A\in \s$ such that $\R(A)\cap {\bR}\neq \emptyset$.  Suppose  $S_n$ is a sequence of densely defined  symmetric operators satisfying  the following conditions.

\noindent (a) $D(A)\subset D(S_n)$.

\noindent (b) There exists a sequence of positive numbers $c_n\ra 0$ such that
\[
\|S_nu\|\leq c_n (\|Au\|+\|u\|),\;\;\forall u\in D(A).
\]
Then  $A+S_n\in \s$ for all $n\gg 0$ and
\[
\rho(A+S_n, A)\ra 0.
\]
\label{prop: rho}
\end{proposition}

\noindent{\bf Proof}\hspace{.3cm}  Set $A_n:=A+S_n$. According to  \cite[Thm.IV.2.24]{Kato} we have
\[
\gamma(A_n, A)\ra 0
\]
 while \cite[Thm. V.4.1]{Kato} implies $A+S_n\in \s$ for  all sufficiently large $n$. Let $\beta\in \R(A)\cap {\bR}$ and  consider a small closed interval $I=[\beta-\ve, \beta +\ve]$ such that $I\subset \R(A)$. Then, using \cite[Thm. VI.5.10]{Kato} we deduce that for  $n$ sufficiently large we have
 \[
 I\subset \R(A_n),\;\;\forall n\gg 0.
 \]
 Pick now a function  $\alpha\in \A$ such that $\alpha(\lambda)\equiv 1$ for $\lambda\geq \beta+\ve$ and $\alpha(\lambda)\equiv 0$ for $\lambda \leq \beta -\ve$. Using \cite[Thm. VI.5.12]{Kato} we deduce
 \[
 \|\alpha(A_n)-\alpha(A)\|\ra 0.
 \]
 We can now invoke  Proposition \ref{prop: riesz} to conclude that $\rho(A_n, A)\ra 0$.  $\Box$

\section{Families of boundary value problems}
\setcounter{equation}{0}
\label{s: 2}

Consider now as in \cite[App. A]{N} the following data.

\bigskip

\noindent $\bullet$ A  compact, oriented Riemannian manifold  $(M, g)$ with boundary  $N=\partial M$  such that  a tubular neighborhood of $N\hra M$ is  {\em isometric} to the cylinder
\[
([0,1]\times N, dt^2 +g_N)
\]
where $g_N$ is a Riemann metric on $N$ and $t$ denotes the outgoing  longitudinal coordinate.

\medskip

\noindent $\bullet$  An Euclidean bundle of  Clifford modules  $E\ra M$ with  Clifford multiplication
\[
{\bc}: T^*M\ra {\rm End}\,(E).
\]
( ${\bc}(\alpha)$ is skew-symmetric for any real $1$-form $\alpha$.) Set $E_0:=E\!\mid_N$

\medskip

\noindent $\bullet$ $D: C^\infty(E)\ra C^\infty(E)$ a  symmetric Dirac operator with principal symbol ${\bf c}$  such that near $N$ it has the form
\[
D= J(\partial_t - D_0),\;\;J:={\bc}(dt)
\]
where $D_0:C^\infty(E_0)\ra C^\infty(E_0)$ is symmetric and independent of $t$.

\medskip

\noindent $\bullet$  A sequence of symmetric endomorphisms of $E$  independent of $t$ near $N$ such that
\[
\|T_n\|_{C^2}\ra 0
\]
and  (near $N$) the endomorphism $JA_n$ is symmetric. Set $D_n:= D+T_n$. Observe that near $N$ $D_n$ has the form
\[
D_n:= J(\partial_t - D_0 -JT_n).
\]
\bigskip

Following \cite{BW}, we consider the family  $\p$ of  admissible  boundary conditions. It  consists  of  zero order, formally  selfadjoint pseudodifferential projectors    with the same principal symbol as the Calderon projector  of $D_0$. The symbol of any $P$ in $P$ commutes with the symbol of $D_0$ so that the commutator $[P, D_0]$ is a zeroth order pseudodifferential operator. We define a metric $\nu$ on $\p$  by setting
\[
\nu(P, Q):= \Bigl\| P-Q\Bigr\|+ \Bigl\|[P-Q, D_0]\Bigr\|
\]
where $\|\bullet\|$ denotes the norm   on the space of  bounded  operators $L^2(E_0)\ra L^2(E_0)$.

Suppose  now that we are given a projector $P\in \p$ and a sequence $(P_n)\subset \p$.  As in \cite{BW}, we can form the  Fredholm selfadjoint operators
\[
A_n: D(A_n)\subset L^2(E)\ra L^2(E), \;\; D(A_n)=\{ u\in H^1(E);\;\;P_nu\!\mid_N=0\}
\]
\[
A_n u = D_n u
\]
and
\[
A: D(A)\subset L^2(E)\ra L^2(E), \;\; D(A)=\{ u\in H^1(E);\;\;Pu\!\mid_N=0\}
\]
\[
A u = Du.
\]

\begin{proposition} If
\begin{equation}
\lim_{n\ra \infty}\nu(P_n, P)=0
\label{eq: proj}
\end{equation}
Then
\[
\lim_{n\ra \infty}\rho(A_n, A)=0.
\]
\label{prop: bvp}
\end{proposition}

\noindent{\bf Proof}\hspace{.3cm} The proof relies on the following   technical result.

\begin{lemma} There exists a sequence  of bounded, invertible operators $U_n: L^2(E)\ra L^2(E)$ such that

\noindent (i) $1-U_n$  and $1-U_n^*$ define  bounded operators  $H^1(E)\ra H^1(E)$

\noindent (ii) $(U_n-1), (U_n-1)^*\ra 0$   in the norm topology on the space of bounded operators $H^s(E)\ra H^s(E)$,  $s=0,1$.

\noindent (iii) $D(A_n)= U_n^*D(A)$, $\forall n$.

\label{lemma: unitary}
\end{lemma}

We will prove  this lemma after we  have finished the proof of Proposition \ref{prop: bvp}.  Set
\[
B_n: =U_n A_n U_n^*.
\]
Observe that $B_n\in \s$ and $D(B_n)=D(A)$. Moreover
\[
\rho(B_n, A_n)=\|\Psi(U_nA_nU_n^*)-\Psi(A_n)\|= \|U_n\Psi(A_n)U_n^*-\Psi(A_n)\|
\]
\[
= \Bigl\|(\,(U_n-1)+1)\Psi(A_n)(\,(U_n-1)+1)^*-\Psi(A_n)\Bigr\| \leq C\|(U_n-1)\|_{L^2, L^2}\cdot\|\Psi(A_n)\|\ra 0
\]
Thus it suffices to show that
\[
\rho(B_n, A)\ra 0.
\]
Observe  that  for all $u\in D(A)$ we have
\[
\|B_nu -Au\|= \|U_n(D+T_n)U_n^* -D\| \leq \| U_n D(U_n^*u -u)\| +\| U_nT_nU_n^*u\|
\]
\[
\leq \|U_n\|_{L^2, L^2}\|D(U_n^*u-u)\|_{L^2} + C\|T_n\|_{C^2}\|u\|_{L^2}\leq C\Bigl(\|(U_n^*-1)u\|_{H^1}+ \|T_n\|_{C^2}\|u\|_{L^2}\Bigr)
\]
\[
\leq C\Bigl(\|(U_n^*-1)\|_{H^1, H^1}\|u\|_{H^1}+ \|T_n\|_{C^2}\|u\|_{L^2}\Bigr)
\]
(use the elliptic estimates in \cite{BW})
\[
\leq C\Bigl\{\|(U_n^*-1)\|_{H^1, H^1}(\|Au\|_{L^2}+\|u\|_{L^2}) + \|T_n\|_{C^2}\|u\|_{L^2}\Bigr\}\leq c_n(\|Au\|+\|u\|)
\]
where $c_n\ra 0$. Thus, the operator $S_n = B_n-A$ satisfies all the  conditions in Proposition \ref{prop: rho}. On the other hand, $A$ has compact resolvent so that $\R(A)\cap {\bR}\neq \emptyset$. We deduce
\[
\rho(A, B_n)=\rho(A, A+S_n) \ra 0. \;\;\Box
\]

\bigskip

\noindent{\bf Proof of Lemma \ref{lemma: unitary}} \hspace{.3cm} Following  the constructions in \cite[I.\S 6.4]{Kato} define
\[
\hat{U}_n:L^2(E_0)\ra L^2(E_0),\;\;\hat{U}_n=P_nP+(1-P_n)(1-P)=2P_nP-(P_n+P)+1
\]
\[
=2(P+R_n)P-(2P+R_n)+1= R_n(2P-1) +1.
\]
$\hat{U}_n$ is a pseudodifferential operator  of order zero with principal symbol $1$. Observe that
\[
 \hat{U}_n^*=PP_n+(1-P)(1-P_n)
 \]
 and, as explained in \cite[I.\S 6.4]{Kato},  $\hat{U}_n^*$ is invertible and maps $\ker{P}$ onto $\ker P_n$.  Observe  moreover that
\begin{equation}
\|\hat{U}_n-1\|_{L^2,L^2}\leq \|R_n\|_{L^2,L^2}\|(2P-1)\|_{L^2, L^2}\ra 0.
\label{eq: hu}
\end{equation}
Next, observe that
\[
[D_0, \hat{U}_n]=[D_0, R_n](2P-1)+2R_n[D_0, P]
\]
 defines a bounded operator $L^2(E_0)\ra L^2(E_0)$ and, using (\ref{eq: proj}) we deduce
\begin{equation}
\bigl\|\, [D_0, \hat{U}_n]\, \bigr\|_{L^2,L^2}\ra 0.
\label{eq: hu1}
\end{equation}
Observe that $\hat{U}_n$ defines in an obvious fashion a bounded operator
\[
\hat{U}_n: L^2(E\mid_{[0,1]\times N})\ra L^2(E\mid_{[0,1]\times N})
\]

 Consider now   a  smooth increasing function
\[
\eta: [0,1]\ra [0,1]
\]
such that $\eta(t)\equiv 0$ for $t<1/4$ and $\eta(t)\equiv 1$ for $t>3/4$.  We can regard $\eta$ as a function on the tubular neighborhood of $N\hra M$ and then  extending it by $0$ we can regard it as a smooth function on $M$. Notice that    if  $u$ is a  section of $E$ then we can regard $\eta u$ as a section of $E\mid_{[0,1]\times N}$.

For any  section of $E$ smooth {\em up to the boundary} define
\[
U_nu = (1-\eta)u + \hat{U}_n(\eta u).
\]
It is clear that $U_nu$ is  smooth up to the boundary.  Notice also  that there exists a constant $C>0$ independent of  $n$ such that
\[
\|U_nu\|_L^2\leq C\|u\|_{L^2}
\]
for any section $u$ smooth up to the boundary. Thus $U_n$ extends to a bounded operator $L^2(E)\ra L^2(E)$. Using (\ref{eq: hu}) we deduce that
\[
\|(U_n-1)\|_{L^2,L^2}\ra 0.
\]
We want to show that  $U_n$ induces a bounded operator  $H^1(E)\ra H^1(E)$ and then estimate the norm of $(U_n-1)$ as a bounded operator $H^1\ra H^1$.

First of all  observe that  the elliptic estimates for $D_0$ imply that there exists  a  positive constant $C$ such that   if  $u$ is smooth up to the boundary then
\[
C^{-1}\|u\|_{H^1([0,1]\times N)}\leq \|\partial_t u\|_{L^2([0,1]\times N)}+ \|D_0u\|_{L^2([0,1]\times N)}\leq C\|u\|_{H^1([0,1]\times N)}
\]
Observe that for any  section  $u$ smooth  up to the boundary we have
\[
\|U_nu - u\|_{H^1(M)} =\| (1-\eta)u + \hat{U}_n (\eta u) -u\|_{H^1(M)}
\]
\[
=\| \hat{U}_n(\eta u) -\eta u\|_{H^1(M)}=\|\hat{U}_n(\eta u)-(\eta u)\|_{H^1([0,1]\times N)}
\]
\begin{equation}
 \begin{array}{c}\leq C\Bigl(\|\hat{U}_n(\eta u) -(\eta u)\|_{L^2([0,1]\times N)} +\|\partial_t\hat{U}_n(\eta u)-\partial_t(\eta u)\|_{L^2([0,1]\times N)} \\
 \\
 + \|D_0\hat{U}_n(\eta u)-D_0(\eta u)\|_{L^2([0,1]\times N)}\Bigr)
 \end{array}
\label{eq: err}
\end{equation}
Using  (\ref{eq: hu}) we deduce
\[
\|\hat{U}_n(\eta u) -(\eta u)\|_{L^2([0,1]\times N)}\leq c_n\|u\|_{L^2(M)},\;\;c_n\ra 0.
\]
To estimate the second term in (\ref{eq: err}) notice first that $[\partial_t , \hat{U}_n]=0$ so that we have
\[
\|\partial_t\hat{U}_n(\eta u)-\partial_t(\eta u)\|_{L^2([0,1]\times N)}= \|\hat{U}_n\partial_t(\eta u)-\partial_t(\eta u)\|_{L^2([0,1]\times N)}
\]
\[
\leq    c_n\|\partial_tu\|_{L^2([0,1]\times N)} \leq c_n\|u\|_{H^1(M)},\;\; c_n\ra 0.
\]
The  estimate  of the third term in (\ref{eq: err}) requires   a bit more work.  Observe that
\[
D_0\hat{U}_n(\eta u)-D_0(\eta u)= [D_0, \hat{U}_n](\eta u) +\hat{U_n}(D_0\eta u)-D_0(\eta u)
\]
\[
=\eta \Bigl([D_0, \hat{U}_n]u + \hat{U}_n(D_0u) - D_0u\Bigr)
\]
so that
\[
\|D_0\hat{U}_n(\eta u)-D_0(\eta u)\|_{L^2([0,1]\times N)}\leq \|\, [D_0, \hat{U}_n]u\|_{L^2([0,1]\times N)}+ \|\hat{U}_n(D_0u) - D_0u\|_{L^2([0,1]\times N)}
\]
(use (\ref{eq: hu}))
\[
\leq c_n(\|u\|_{L^2([0,1]\times N)} + \|D_0u\|_{L^2([0,1]\times N)}) \leq c'_n\|u\|_{H^1(M)},\;\;c'_n\ra 0.
\]
We have thus    found a sequence of positive numbers $c_n\ra 0$ such that
\[
\|U_n u-u\|_{H^1(M)}\leq c_n\|u\|_{H^1(M)}
\]
for every section $u$ smooth up to the boundary. This shows that $U_n$ induces a bounded operator $H^1(M)\ra H^1(M)$ and  moreover,
\[
\|U_n- 1\|_{H^1, H^1}\leq c_n \ra 0.
\]
One can prove a similar statement   concerning $U_n^*$.   Clearly $U_n$ is invertible being so close to $1$.   Since $\ker P_n=\hat{U}_n^*(\ker P)$ we deduce that $D(A_n)=U_n^* D(A)$. Lemma \ref{lemma: unitary} is proved. $\Box$

\section{Classifying spaces for $K$-theory}
\setcounter{equation}{0}
\label{s: 3}

For clarity purposes    we will consider only a special case, that of the functor $KO^1$. To  discuss  the other functors $KO^n$ one should use the bigraded  Karoubi functors $KO^{p,q}$ as we did in \cite{N}. The proof is only notationally more complicate.

  Denote by $\f\subset \s$ (resp. $\mathcal{BF}\subset \mathcal{BS}$, $[\mathcal{BF}]\subset[\mathcal{BS}]$) the subspace of selfadjoint Fredholm operators. $[\mathcal{BF}]$ has three connected components. Two of them $[\mathcal{BF}_\pm]$, are contractible while the third, $[\mathcal{BF}_0]$  is a classifying space for $KO^1$ (see \cite{AS, BW,  Kar}). We deduce that $(\f, \rho)$ consists of three components
\[
\f_\pm:=\Psi^{-1}([\mathcal{BF}_\pm]),\;\;\f_0:=\Psi^{-1}([\mathcal{BF}_0])
\]
and $(\f_0, \rho)$ is a classifying space for $KO^1$.

Observe that $H\oplus H$ is a symplectic space with complex structure
\[
J=\left[
\begin{array}{cc}
0 & -1_H\\
1_H & 0
\end{array}
\right]
\]
and $\Lambda_0:=H\oplus 0$ is a Lagrangian subspace. Define $\mathcal{FL}_0$ the set of Lagrangian subspaces $\Lambda\subset H\oplus H$ such that $(\Lambda_0, \Lambda)$ is a Fredholm pair. We topologize $\mathcal{FL}_0$ using the gap distance $\delta$. The space $(\mathcal{FL}_0, \delta)$  is also a classifying space for $KO^1$ (see \cite{N0}).

There is a natural $1-1$ map
\[
\Gamma:\f_0 \ra \mathcal{FL}_0,\;\;A \mapsto \Gamma_A.
\]
According to Lemma \ref{lemma: ga-ri} the  map $\Gamma: (\f_0, \rho)\ra (\mathcal{FL}_0, \delta)$ is continuous.

\begin{theorem} The map
\[
\Gamma: (\f_0, \rho)\ra (\mathcal{FL}_0, \delta)
\]
is a weak homotopy equivalence.
\label{th: cl}
\end{theorem}

\noindent{\bf Proof}\hspace{.3cm} Fix $A_0\in \f_0$. We have to show that for every $n>0$ the induced map
\[
\Gamma_*:\pi_n(\f_0, A_0)\ra \pi_n(\mathcal{FL}_0, \Gamma_{A_0})
\]
is an isomorphism.  Observe first that, according to Bott periodicity,
\[
\pi_n( \mathcal{FL}_0, \Gamma_{A_0})\in \mathcal{G}:=\Bigl\{0, {\bZ}, {\bZ}_2\Bigr\}.
\]
The groups in the family $\mathcal{G}$ have a remarkable property. If $G\in \mathcal{G}$ and $\vfi: G\ra G$ is a surjective morphism then $\vfi$ is an isomorphism.

In \cite[\S 5.3]{N}, using the symplectic reduction morphism it is shown that the morphism $\Gamma_*$ is surjective provided the (general)  Floer families are  $\rho$-continuous.  This continuity was established  in  Proposition \ref{prop: bvp}. Theorem \ref{th: cl} is proved. $\Box$

\begin{remark}  In \cite{N} we claimed that the map $\Gamma: (\f_0, \gamma)\ra (\mathcal{FL}_0, \delta)$ is a weak homotopy equivalence when in fact the  arguments there, detailed in this paper, prove this only for the stronger $\rho$-topology. This has no effect on the results  of  \cite{N} but one    {\ae}sthetical  question still lingers. Is the space $\mathcal{F}_0$ equipped with the gap  topology a classifying space for $KO^1$? If the answer is yes (which we continue to belive to be the case) then our claim in \cite{N} is true.
\end{remark}

\end{document}